\documentclass{amsart}
\usepackage{amsfonts}

\theoremstyle{plain}

\numberwithin{equation}{section}
\usepackage{graphicx} 

\DeclareGraphicsExtensions{.pdf, .jpg}

\begin{document}

\title{$q$-exponential behavior of expected  aggregated supply curves in deregulated electricity markets}
\author{C. Cadavid}
\address{Carlos Cadavid, {\small {Universidad EAFIT, Departamento de
Ciencias Matem\'aticas, Bloque 38, Office 417 }}Carrera 49 No. 7 Sur -50,
Medell\'{\i}n, Colombia (57)(4)-2619500.}
\email{ccadavid@eafit.edu.co.}
\author{M. E.  Puerta}
\address{Mar\'ia E. Puerta, {\small {Universidad EAFIT, Departamento de
Ciencias Matem\'aticas, Bloque 38, Office 419 }}Carrera 49 No. 7 Sur -50,
Medell\'{\i}n, Colombia (57)(4)-2619500.}
\email{mpuerta@eafit.edu.co}
\author{J. D. V\'{e}lez}
\address{Juan D. V\'{e}lez, {\small{ Universidad Nacional de Colombia, Departamento de Matem\'aticas, Medell\'in, Colombia}}}
\email{jdvelez@unal.edu.co}
\author{Juan F. Garc\'ia}
\address{Juan Fenando Garc\'ia Pulgar\'in, {\small {Universidad EAFIT, Departamento de
Ciencias Matem\'aticas, Bloque 38, Office 417 }}Carrera 49 No. 7 Sur -50,
Medell\'{\i}n, Colombia (57)(4)-2619500.}
\email{jgarci77@eafit.edu.co}
\keywords{}
\thanks{}

\maketitle

\begin{abstract}
It has been observed that the expected aggregated supply curves of the colombian electricity market present  $q$-exponential behavior.  
The purpose of this article is to present evidence supporting the fact that $q$-exponentiality is already present in the expected aggregated supply 
curve of certain extremely simple idealized deregulated electricity markets where illegal interaction among competing firms is precluded.
 
\end{abstract}

\section{Introduction}

In the last two decades, more and more countries have adopted a deregulated electricity market policy. In these countries, an inverse 
auction mechanism is implemented each day in order to determine, for the next day, which generating units will be operating, how much electricity 
will be supplied by each one of them, and its unitary price. Since the electrical sector is of utmost importance in any country's economy, it is important 
to understand various aspects of the corresponding market. One of such aspects is the \emph{expected aggregated supply curve}, which will be refered to as \emph{supply curve} for short,  resulting from the competition among generating firms. 
It has been observed that the supply curve for each of the years in the period $2005$-$2011$, in the colombian electricity market (a deregulated one), 
are well fitted by suitably scaled $q$-exponential curves, where the values of $q$ vary from year to year and are generally different from $1$ \cite{gabriel}. 
The presence of $q$-exponential behavior or $q$-deformed behavior in general, has been observed in quantities arising in similar contexts. In \cite{kristo} certain quantities 
associated to the Czech Republic public procurement market (an auction type market), for the period $6/2006$-$8/2011$, are analyzed. Specifically, those authors look at the behavior of three quantities, 
namely, the probability that a contract has $x-1$ bidders or more, the probability that a  public procurement selling agent makes $x$ euros or more during the observed period ($6/2006$-$8/2011$), and the 
probability that a contracting authority spends $x$ or more euros during the observed period, as functions of $x-1$ in the first case, and $x$ in the second and third cases. 
They observe that the first function is very well approximated by certain 
exponential, while the second and third functions, are very well approximated by certain power functions. 
Now, exponential functions and power functions correspond to $q$-exponentials, with $q=1$ in the first case and $q\neq 1$ in the second one. 
In \cite{chininis} the probability distribution of certain quantities associated to the opening call auctions in the chinese stock market are studied. In particular, it is shown that the probability distribution of order sizes is well fitted by
certain $q$-Gamma function. \\

Even though the presence of $q$-deformed behavior, with $q\neq 1$, may be taken as evidence of complexity in the dynamics governing the corresponding processes, illegal interaction between participants and so on, it is still valid to ask whether there is already 
$q$-deformed behavior, with $q\neq 1$, in situations having low complexity and complete transparency. 

 In this article we study the presence of $q$-deformed behavior in the seemingly most simple situation arising in deregulated electricity markets. Concretely, we consider an electricity market composed of a number of generating firms, each one having a single generating set capable of supplying (after normalization) one electricity unit per day whose production cost is assumed (after normalization) to be zero. In addition, the demand is a discrete random variable whose probability distribution is known by all firms.  As was mentioned before, the firms compete each day in a reverse auction in order to be selected as an electricity supplier for the next day. This simple market can be thought as a repeated game in the framework of classical game theory, i.e. one in which players are rational and have complete knowledge of the game structure. Therefore, we assume that players place their bids according to a Nash equilibrium. After calculating and solving a differential equation whose solution is the unique symmetric Nash equilibrium possesed by this game,  we estimate the supply curve by simulating repetitions of this game. We observed that for various background probability distributions for the demand, and a range of numbers of generating firms, the supply curves are well fitted by $q$-exponential curves, usually with $q\neq 1$, and that for fixed background distribution, the value of $q$ gets closer to $1$ as the number of generating firms increases.  Towards the end of the paper we propose a plausible explanation for these phenomena.

\section{General model}

In this section we closely follow \cite{von der Fehr}. We assume that the electricity generating system is formed by $N$ electricity generating firms $g_1,\ldots,g_N$, and that each firm $g_n$ has $m_n$ generating units or sets  $s_{n1},\ldots,s_{nm_n}$. Let $M=\sum_{n=1}^N m_n$ be the total number of sets of the whole generation system. Each set $s_{ni}$ can produce at most $k_{ni}$ electrity units per day and has a cost function $c_{ni}:[0,\infty)\to \mathbb{R}$, assigning to each $q\in [0,\infty)$, the cost $c_{ni}(q)$ of producing $q$ electricity units using $s_{ni}$. Let $K=\sum_{ni}k_{ni}$ be the total daily capacity of the generating system. We assume that the amount of electricity demanded by the consumers in one day is a random variable $d$ that distributes according to certain cumulative distribution function $G$ supported on some interval $[\underline{d},\overline{d}]\subset [0,K]$. Policies for the market demand that firms can only offer unitary prices lying within certain interval $[\underline{p},\overline{p}]$. It is assumed that all the information above is common knowledge. The mechanism by which electricity is bought is the following.  Each day (day $t$), each $g_n$ secretely submits a vector $(p_{n1},\ldots,p_{nm_n})$ to a coordinating entity, expressing its willingness to produce the next day (day $t+1$), using set $s_{ni}$, any amount of electricity in the interval $[0,k_{ni}]$, and to sell it at a unitary price of $p_{ni}$. At this point we emphasize the fact that the generators decide on their price vectors knowing the distribution $G$ but without knowing the particular value it takes the next day. Once the coordinating entity receives all these $N$ vectors, it chooses a \emph{ranking} of the prices, i.e. a bijective function $r:\{(n,i):n=1,\ldots,N, i=1,\ldots,m_n\} \to \{1,\ldots,M\}$ such that $r(n,i)<r(m,j)$ whenever $p_{ni}<p_{mj}$. Notice that this choice is not unique whenever there exist pairs $(n,i)\neq (m,j)$ such that $p_{ni}=p_{mj}$. More precisely, if $\{A_1,\ldots,A_l\}$ is the partition induced on the set $\{(n,i):n=1,\ldots,N, i=1,\ldots,m_n\}$ by the equivalence relation $(n,i)\equiv (m,j)$ iff $p_{ni}=p_{mj}$, the number of choices is $\vert A_1 \vert ! \cdots \vert A_l \vert !$. Here $\vert A_j \vert$ denotes the number of elements of the set $A_j$. This number will be denoted by $R(p_{11},\ldots,p_{Nm_N})$. The coordinator randomly selects one of the possible rankings according to a uniform distribution, i.e. with probability $\frac{1}{\vert A_1 \vert ! \cdots \vert A_l \vert !}$.  Once a ranking $r$ is chosen, the coordinator looks at the actual demand of electricity $d$ for day $t+1$. It is useful at this point to rename the generating units and all the entities associated to them, according to their rank, i.e.  to denote the unit  $s_{ni}$ by $s_{r(n,i)}$, $k_{ni}$ by $k_{r(n,i)}$, and so on.  The coordinator considers the numbers $K_0=0$ and $K_j=\sum_{a=1}^jk_{a}$ for $j=1,\ldots,M$, and determines $\rho=\max \{j:K_{j-1}<d\}$. Then the coordinator dictates that on day $t+1$, each unit $s_j$ with $j<\rho$ will supply its full capacity $k_j$, that unit $s_{\rho}$ will supply $d-K_{\rho-1}$ units of electricity, and that each unit of electricity will be paid at price $p_{\rho}$. It is understood that units $s_j$ with $j>\rho$ will not supply any electricity and therefore will not receive any payment. The utility for firm $g_n$ is therefore given by 
\begin{equation}
u_n=\delta_n(\rho)\left[(d-K_{\rho-1})p_{\rho}-c_{\rho}(d-K_{\rho-1}) \right]  +  \sum_{j=1}^{\rho-1}\delta_n(j)\left[k_jp_{\rho}-c_{j}(k_j)\right]
\end{equation}
where $\delta_n(l)$ is $1$ if $s_l$ belongs to firm $g_n$ and is zero otherwise.

Since we will consider the case when (i) each generator $g_n$ has a single generating unit $s_{n1}$, (ii) each generating unit has capacity $1$, (iii) all cost functions are $c_{n1}=0$, (iv) the demand is a discrete random variable taking values $\{1,\ldots,N\}$ , and (v) $\underline{p}=0$; \emph{these five conditions will be assumed to hold for the rest of the paper}. In this case it is important to slightly modify our notation.  Instead of writing $s_{n1},k_{n1},c_{n1},p_{n1}$ we will write $s^n,k^n,c^n,p^n$, respectively. The renaming of units and associated entities induced by a ranking $r$ keeps being the same as before, i.e.  $s^n, k^n, c^n,p^n$ are renamed as $s_{r(n,1)}, k_{r(n,1)}, c_{r(n,1)}, p_{r(n,1)}$, respectively.  

The assumption that firms want to maximize their utilities, naturally leads to the interpretation of the whole situation as the repetition of a three stage game: first, firms choose their prices in the interval $[0,\overline{p}]$, then nature chooses a value for the demand according to some probability distribution $\pi_i=Pr(d=i)$, $i=1,\ldots,N$, and finally the coordinator makes the dispatchment according to the rules already explained.  Notice that the latter step involves a random choice in order to resolve ties. The utilities for each firm depend on the prices being offered by all firms, the value taken by the demand and the particular ranking chosen by the coordinator.  This game can thought as a random experiment, whose sample space is
\begin{equation} 
\Omega=\left\{ (p_1,\ldots,p_N,d,r): \text{each} \ p_n\in [0,\overline{p}], d\in \{1,\ldots,N\}, r \ \text{a ranking for} \ (p_1,\ldots,p_N)\right\}
\end{equation}
with probability density function $f$ (pdf) having the form 
\begin{equation}
f(p_1,\ldots,p_N,d,r)=f_1(p_1)\ldots f_N(p_N).\pi_d.\frac{1}{R(p_1,\ldots,p_N)}\,.
\end{equation}
Here $f_1,\ldots,f_N$ are some pdf's. In this setting the utility functions $u_1,\ldots,u_N$ become random variables. The functions $f_1,\ldots,f_N$ according to which firms choose their prices, form a Nash equilibrium, if for each $n$, the expected utility of firm $n$, 
\begin{equation}
E_{f_1,\ldots,f_n,\ldots,f_N}(u_n)\geq E_{f_1, \ldots,\tilde{f_n},\ldots,f_N}(u_n) 
\end{equation}
for any other pdf $\tilde{f}_n$.
It is equivalent but much easier to work with the corresponding cummulative distribution functions $F_1,\ldots,F_N$. These will be nondecreasing functions defined on $[0,\overline{p}]$, having value $0$ at $p=0$ and $1$ at $p=\overline{p}$, and admitting (jump) discontinuities.  

According to proposition 6 in \cite{von der Fehr} $F_1=\ldots =F_N=F$ is a (symmetric) Nash equilibrium if and  only if $F$ is a nondrecreasing function defined on a closed interval  $[p^m,\overline{p}]$ with $0<p^m<\overline{p}$, such that $F(p^m)=0$, $F(\overline{p})=1$, and satisfying the differential equation 

\begin{equation}\label{derivada}\sum_{i=1}^{N}\pi_i\left\{ H_i'p + H_i + p G_i'(p)\right\} =0\end{equation} 

where

\begin{equation*} H_i=\left( \begin{array}{ccc} N-1 \\ i-1 \end{array} \right) F(p)^{i-1}(1-F(p))^{N-i}\end{equation*}

and 
$$G_i(p)=\sum_{k=i-1}^{N-1}\left( \begin{array}{ccc} N-1 \\k\end{array} \right)F(p)^k(1-F(p))^{N-1-k}\,.$$

We remark that this differential equation derives from the general fact that a profile of mixed strategies $F_1,\ldots,F_N$ is a Nash equilibrium if and only if for each $n$,  the (expected) profit of firm $g_n$  assuming that it plays the pure strategy $p_n=p$ while any other firm $g_{n'}$ plays according to $F_{n'}$, is independent of $p$. We observe that this differential equation can be solved by separation of variables.

\section{Experiments and results}

In order to survey the behavior of supply curves in general, we wrote routines in Mathematica 10.1 and used them for testing various scenarios.    
The main routine, called GranProgramaP,  assumes $\bar{p}=100$ and for given distribution of background demand, and positive integers $N$ and $T$, it first determines the differential equation \ref{derivada}, finds the solution $F$ to this differential equation subject to the condition $F(100)=1$, and then, for each $1\leq i \leq T$, it produces an $N$-vector $v_i$, by first making $N$ independent (random) choices of numbers in the interval $[0,100]$ according to $F$, sorting these numbers in increasing order, and finally computing the average $N$-vector $V=\frac{1}{T}\sum_{i=1}^{T} v_i$. 
Finally, the Mathematica command NMinimize is applied to find, among the curves $y=f(x)=\alpha \exp_q(\beta x)$ with $\alpha,\beta$, $q>0$, and $N-\frac{1}{2}<\frac{1}{\beta (q-1)}$ whenever $q>1$, one that (approximately) minimizes the quadratic error
$$\sum_{j=1}^{N} \left(V_j-f\left(j-\frac{1}{2}\right)\right)^2\,.$$

Next we present the considered scenarios and the results obtained in each case.  The reason for the kind of experiment performed in each of these scenarios is clarified in section \ref{analysis}. The following four experiments have the following common structure. 
We assume a market with $N$ generating firms, each of them capable of producing $1/N$ units per day at production cost zero, and with demand being a discrete random variable taking values $i/N$ for $i=1,\ldots,N$, with probability $$\pi_i=Pr(\text{demand}=i/N)=  \int_{\frac{i-1}{N}}^{\frac{i}{N}} \! h(x) \, \mathrm{d}x,$$ for a fixed probability distribution $h$ supported on the interval $[0,1]$. We refer to $h$ as the \emph{background distribution} for the demand.
This is of course equivalent to a market with $N$ generating firms, each of them capable of producing $1$ unit per day at production cost zero, and with demand being a discrete random variable taking values $i=1,\ldots,N$, with probability 
$$\pi_i=Pr(\text{demand}=i)=  \int_{\frac{i-1}{N}}^{\frac{i}{N}} \! h(x)  \, \mathrm{d}x.$$ So the experiments only differ by the choice of background distribution.  In all four experiments we set $T=10000$ and for NMinize we selected the \emph{RandomSearch} method with $500$ search points.

\subsection{Experiment 1}

In this experiment we took $h(x)=1, \  0\leq x \leq 1$ as background distribution for the demand.
 Table \ref{tabla1}  shows the results obtained for values of $N$ from $5$ to $28$. 

\begin{table}[h]
\begin{tabular}{ |p{1cm}||p{2cm}|p{3cm}|p{2cm}|p{3cm}|} 
 \hline
 \multicolumn{5}{|c|}{Uniform demand} \\
 \hline
 $N$     & $q$ &  $\alpha$ & $\beta$ & quadratic error \\
 \hline
5 & 0.863233 & 2.9002 & 0.823226 & 0.0193731 \\
6 &  0.89449 &  1.2304 & 0.859379 &  0.0672908 \\
7 & 0.892939 & 0.437755 & 0.980874 & 0.00743362 \\ 
8 &  0.916562 &  0.183676 & 0.972093  & 0.00345962 \\
9 &  0.933164 & 0.0836263 & 0.955353 &  0.00345632 \\
10 & 0.93709 & 0.0322543 & 0.997464 &  0.00232163 \\ 
11 & 0.939185 &  0.0121861  & 1.04244  &  0.00379461\\ 
12 & 0.947345 & 0.00562502 & 1.02362 & 0.00259078  \\
13 & 0.943343  & 0.00147738 & 1.1462 &  0.0115656  \\ 
14 & 0.957351 & 0.000978965 & 1.02787 &  0.00465776\\  
15 &  0.962617 & 0.000435722 & 1.01207 &  0.0013087 \\ 
16 & 0.987828 & 0.000814185  & 0.765452 &  0.0169305 \\
17 & 0.969912  &   0.000107791  &  0.972005 & 0.00665632\\ 
18 & 0.950759 &  $4.94466\times10^{-6}$ & 1.4119 &  0.0184282\\ 
19 & 0.968327 & 0.0000111811 & 1.06863 &  0.00827299\\ 
20 &  0.973942  &  $8.62225\times10^{-6}$  &  0.988244 &  0.00439083\\ 
21 &  0.991221 & 0.0000229271  & 0.76225 & 0.0286424 \\ 
22 & 0.982056 & $2.77873\times10^{-6}$ & 0.908473 &  0.0226434\\ 
23  & 0.984581 &  $1.54111\times10^{-6}$  &  0.882888 &   0.0220357\\ 
24 &  0.992261  & $2.3351\times10^{-6}$   & 0.769326 & 0.0998453\\
25 & 0.995941    & $2.15478\times10^{-6}$  & 0.717738 &  0.0957674 \\ 
26 & 0.993057  & $6.32775\times10^{-7}$ & 0.761542 &   0.087408\\   
27 & 1.0081 & $2.26665\times10^{-6}$ & 0.596923 &  0.154767\\  
28 &  1.01065 &  $1.5224\times10^{-6}$ & 0.575789 &  0.199678\\ 
\hline
\end{tabular}
\caption{}
\label{tabla1}
\end{table}

\subsection{Experiment 2}

In this experiment we took $h(x)=3x^2, \ 0\leq x \leq 1$ as background distribution for the demand.
 Table \ref{tabla2}  shows the results obtained for values of $N$ from $5$ to $28$.

\begin{table}[h]
\begin{tabular}{ |p{1cm}||p{2cm}|p{3cm}|p{2cm}|p{2.5cm}| } 
 \hline
 \multicolumn{5}{|c|}{Demand's f.d.p based on $3x^2$, $0\leq x \leq 1$} \\
 \hline
 $N$     & $q$ &  $\alpha$ & $\beta$ & quadratic error \\
 \hline
5 & 0.786277 & 9.24399 & 0.562637 &  0.000517729\\
6 & 0.848111 & 4.08727 & 0.632552 &  0.0713172\\  
7 & 0.882261 & 1.72287 & 0.691693 &  0.0499626 \\
8 & 0.905359 &  0.724865 & 0.738601 & 0.012073\\  
9 & 0.910659  & 0.289822 & 0.805672 &  0.00154287\\ 
10 & 0.923683 &  0.114552 &  0.843706 & 0.0121201\\ 
11 & 0.940745 & 0.0609441 & 0.80781 &  0.0103708\\ 
12 & 0.952511 & 0.0259056 & 0.810255 &  0.00760159\\
13 & 0.953217 & 0.0102512 & 0.850152 &  0.00430776\\ 
14 & 0.962814 & 0.00485199 & 0.829615 & 0.0294654\\ 
15 & 0.959879 &  0.00173207 & 0.888678 & 0.0151036\\ 
16 &  0.956018 & 0.000481014 & 0.983035 & 0.00436682\\ 
17 & 0.962618 & 0.000285168 & 0.934857 &  0.0044448\\ 
18 &   0.967372 & 0.000148451 & 0.910055 & 0.0162562 \\ 
19 & 0.974735 &  0.0000893455 &  0.85638 & 0.00924199\\ 
20 & 0.973182  & 0.0000324616  & 0.896141 &  0.0112569\\  
21 & 0.962477 & $2.70374\times10^{-6}$  & 1.14381 & 0.0111554\\ 
22 &  0.990649 & 0.0000259162 & 0.725679 &  0.0354996\\   
23 &  0.975639 & $2.1441\times10^{-6}$ & 0.938534 &   0.0101328\\   
24 & 0.975775  & $8.8424\times10^{-7}$ & 0.955251 &  0.0344693\\ 
25 & 0.983022 &  $9.65909\times10^{-7}$ & 0.850342 &   0.0896586\\   
26 &   0.997108 & $3.26761\times10^{-6}$  & 0.667209 &  0.0643818\\ 
27 & 1.00178 & $2.30702\times10^{-6}$ & 0.63024 &  0.134193\\
28 & 0.991328 & $3.02261\times10^{-7}$ & 0.750202 & 0.0612747\\ 
\hline
\end{tabular}
\caption{}
\label{tabla2}
\end{table}

\subsection{Experiment 3} 

In this experiment we took $h(x)=3(x-1)^2, \ 0\leq x \leq 1$ as background distribution for the demand.
Table \ref{tabla2}  shows the results obtained for values of $N$ from $5$ to $28$. 
\begin{table}[h]
\begin{tabular}{ |p{1cm}||p{2.5cm}|p{3 cm}|p{3 cm}|p{3cm}| } 
 \hline
 \multicolumn{5}{|c|}{Demand's f.d.p based on $3(x-1)^2$, $0\leq x \leq 1$} \\
 \hline
 $N$     & $q$ &  $\alpha$ & $\beta$ & quadratic error \\
 \hline
5 &  1.09464 & 0.206346 & 0.848798 & 0.0024282\\ 
6 & 1.09723 &  0.0711591  & 0.798572 &  0.00177619\\  
7 & 1.10279 & 0.0319655 & 0.730883 & 0.00488975\\  
8 &  1.11106 &  0.0156945 &   0.660128 & 0.00499712\\  
9 & 1.09178 & 0.00541504 &  0.679568 & 0.00123698\\ 
10 & 1.09736 & 0.00281438 & 0.625793 & 0.000942782\\ 
11 & 1.09776 & 0.00163252  & 0.58577 & 0.000979379\\   
12 & 1.06551 &  0.000271507 & 0.690339 & 0.000269798\\ 
13 & 1.0952 & 0.000422896 & 0.538375 &  0.00126523\\ 
14 & 1.08674 & 0.000198589  &  0.538884  &  0.000486154\\   
15 & 1.09073 & 0.00012446 & 0.502844 &  0.00021681\\ 
16 & 1.06933 & 0.0000260694 &  0.565947 & 0.000614531\\ 
17 & 1.08805 & 0.0000376827 & 0.471669 & 0.000287819\\ 
18 & 1.08563 & 0.0000213679  & 0.460938 &  0.000232498\\
19 & 1.09183 & 0.0000198484 & 0.422582 &  0.00157232\\ 
20 & 1.10914 & 0.0000282786  & 0.36162 &  0.00176623\\ 
21 & 1.0921 & $8.29756\times10^{-6}$  &  0.392186  &  0.000280447\\ 
22 & 1.09542 & $6.87487\times10^{-6}$  & 0.368959  & 0.000318908\\ 
23 & 1.09237 & $4.13338\times10^{-6}$ & 0.364318 & 0.00249175\\ 
24 & 1.08245 & $1.33427\times10^{-6}$  &  0.383466 &   0.000381974\\ 
25 &  1.09711 & $2.89082\times10^{-6}$ & 0.32872  &  0.00305721\\  
26 & 1.08742 & $1.0406\times10^{-6}$   & 0.344124 & 0.00312855\\  
27 & 1.11564 & $5.33944\times10^{-6}$ & 0.269265 & 0.00625656\\   
28 & 1.11006 & $2.83686\times10^{-6}$  & 0.271704 &  0.00357036\\
\hline
\end{tabular}
\caption{}
\label{tabla3}
\end{table}

\subsection{Experiment 4} 
\begin{figure}
\includegraphics[bb=20 20 300 200,scale=0.5]{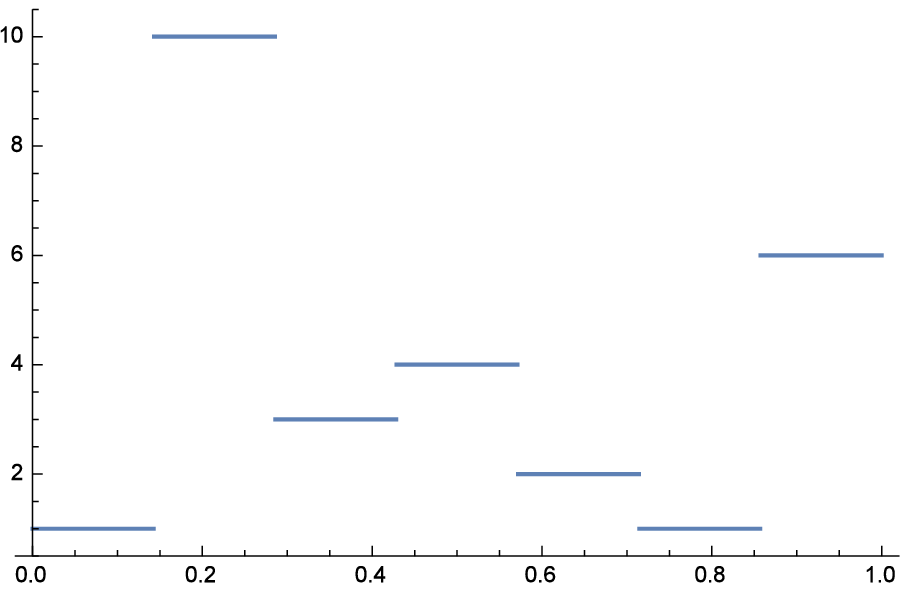}
\caption{}
\label{demandad}
\end{figure}

So far we have considered background demands distributed according to power-like functions. It is important to experiment with more general background distributions for the demand. 
Let us consider for instance a background demand distributed as some randomly selected piecewise constant function $g$ whose graph is depicted in Figure \ref{demandad}. 

 Table \ref{tabla4}  shows the results obtained for values of $N$ from $5$ to $28$.

\begin{table}[h]
\begin{tabular}{ |p{1cm}||p{2cm}|p{2 cm}|p{2.5 cm}|p{2.5cm}| } 
 \hline
 \multicolumn{5}{|c|}{Demand's f.d.p based on $g(x)$, $0\leq x \leq 1$} \\
 \hline
 $N$     & $q$ &  $\alpha$ & $\beta$ & quadratic error \\
 \hline

5 & 0.72133 & 2.954427 &  1.129140 & 0.4091208\\  
6 & 0.67358 & 0.821282 & 1.861038 & 1.8020197\\ 
7 & 0.64121 & 0.038771 & 6.119756 &  2.8822813\\   
8 & 0.69959 & 0.005481 & 7.269161 & 1.2973836\\
9 & 0.73581 & 0.000718 & 8.862518 &  1.39342351\\ 
10 & 0.77651 & 0.000579 & 5.932255 &  0.6524641\\ 
11 &  0.80561 & 0.000245 & 5.070716 &  0.4075098\\ 
12 & 0.83691 & 0.000347 & 3.309873 &  0.0719115\\  
13 & 0.85802  & 0.000301  &  2.624487 &  0.0562562\\  
14 &  0.87835  & 0.000241  &  2.143395 &  0.1501309\\ 
15 & 0.90866  & 0.000540  &  1.418274  &    0.2352072\\ 
16 &  0.92360 & 0.000395 & 1.247219  &  0.2155978\\  
17 &  0.94694 & 0.000562  &  0.961471  &  0.26144589\\ 
18 & 0.97006 & 0.000738  &  0.762271  & 0.2988319\\  
19  & 0.98948 & 0.000945 &  0.628611 &   0.2476883\\  
20 & 0.99752 & 0.000670 & 0.587217 &  0.2163773\\ 
21 & 1.01654 &  0.000856  & 0.491904  &  0.1813817\\  
22  & 1.02361 & 0.000688  &  0.460020 & 0.2267202\\ 
23 & 1.02211 & 0.000373  & 0.463638 & 0.2033697\\  
24 & 1.02958 & 0.000297 &  0.432238 & 0.1293418\\  
25  & 1.02354 &  0.000119 & 0.457114 &  0.0546924\\  
26  & 1.03403 & 0.000130   & 0.409306 &  0.0751101\\  
27  & 1.04037 & 0.000121  & 0.381030 &  0.0664762\\ 
28  &  1.03368 & 0.000044 & 0.405377 & 0.0281123\\ 
\hline
\end{tabular}
\caption{}
\label{tabla4}
\end{table}

\section{Analysis of results and conclusions} \label{analysis}

\begin{figure}
\includegraphics[bb=20 20 300 200,scale=0.5]{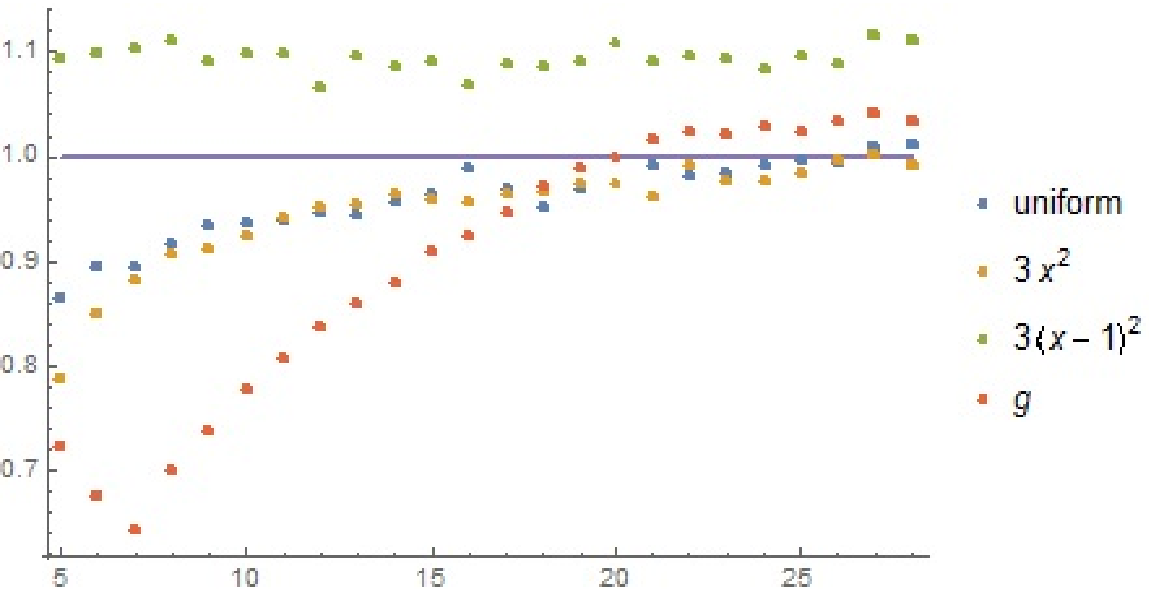}
\caption{}
\label{behaviorofq}
\end{figure}

The results obtained in the previous section strongly indicate the presence of $q$-exponential behavior, with $q\neq 1$, of the supply curve in the various tested scenarios. The quadratic errors oscillated between $0.0013087$ and $0.199687$ in experiment 1, between $0.000517729$ and $0.134193$ in experiment 2, between $0.00021681$ and $0.00625656$ in experiment 3, and between $0.0281123$ and $2.882283$ in experiment 4. These errors are very small in comparison with the quantities being approximated.
It is important to observe that the emergence of $q$-exponential behavior in these experiments cannot be due to any illegal behavior, since the mathematical model assumes that prices are chosen independently by the generating firms. Figure \ref{behaviorofq}

insinuates an interesting tendency: for fixed background demand, the value of $q$ gets closer to $1$ as $N$ gets larger. These observations suggest the following interpretation. When the number of competing firms $N$ is small, we can consider each one of them as a maximally strong coalition of a large number of minifirms. This allows us to consider the small number of competing firms as a huge number of competing firms divided into $N$ groups so that firms belonging to one group strongly interact mutually and firms belonging to different groups do not interact. Therefore $N$ getting larger has the effect of breaking coalitions, reducing in this way the level of interaction, thus making $q$ closer to $1$. We infer that as $N$ gets larger, the market gets closer to a perfect competition market. This suggests the possibility of modeling these markets as statistical mechanical systems having nontrivial particle interaction \cite{tsallis}, and justifying in this way the presence of $q$-exponential behavior in this setting.

\end{document}